\documentclass[10pt]{amsart}
\usepackage{amsxtra,latexsym,amssymb}
\usepackage{epsfig,rotate,amsthm}

\def\qed{{\hfill $\Box$}}

\def\N{{\mathbb N}}
\def\Z{{\mathbb Z}}
\def\C{{\mathbb C}}

\theoremstyle{theorem}
\newtheorem{thm}{Theorem}[section]
\newtheorem{cor}{Corollary}[section]
\newtheorem{prop}{Proposition}[section]

\newtheorem{lem}{Lemma}[section]
\theoremstyle{definition}
\newtheorem{defn}{Definition}[section]
\theoremstyle{remark}

\newtheorem{rem}{\bf Remark}[section]

\begin{document}

\begin{center}
{\Large {\bf ON REPRESENTATIONS OF QUANTUM GROUPS
    $U_{q}(f_{m}(K,H))$}}\\
\vskip .1in
{\footnotesize \it{To Prof. Yingbo Zhang on the occasion of her 60th birthday}}
\end{center}
\vskip .2in
\begin{center}
Xin Tang $^{a}$ and Yunge Xu $^{b,}\footnote{Yunge Xu is partially supported by 
NSFC under grant 10501010.}$

\bigskip

{\footnotesize a. Department of Mathematics $\&$ Computer Science, Fayetteville State University,\\
 Fayetteville, NC 28301, U.S.A \hspace{5mm} E-mail:
xtang@uncfsu.edu

b. Faculty of Mathematics $\&$ Computer Science,
Hubei University,\\ Wuhan 430062, P.R.China \hspace{5mm} E-mail:
xuy@hubu.edu.cn}

\end{center}

\bigskip

\begin{center}

\begin{minipage}{12cm}
{\bf Abstract}: 
We construct families of irreducible representations for a class 
of quantum groups $U_{q}(f_{m}(K,H)$. First, we realize these quantum 
groups as Hyperbolic algebras. Such a realization yields natural 
families of irreducible weight representations for $U_{q}(f_{m}(K,H))$. 
Second, we study the relationship between $U_{q}(f_{m}(K,H))$ and
$U_{q}(f_{m}(K))$. As a result, any finite dimensional weight
representation of $U_{q}(f_{m}(K,H))$ is proved to be completely
reducible. Finally, we study the Whittaker model for the 
center of $U_{q}(f_{m}(K,H))$, and a classification of all 
irreducible Whittaker representations of $U_{q}(f_{m}(K,H))$ 
is obtained.
\medskip


\medskip


\end{minipage}

\end{center}
                 
\section*{0. Introduction}
Several generalizations (or deformations) of the quantized enveloping 
algebra $U_{q}(sl_{2})$ have been extensively studied in \cite{BW, HZ, JZ, JWS, JWZ}. Especially in \cite{JWZ}, a general class of algebras $U_{q}(f(K))$ 
(similar to $U_{q}(sl_{2})$) was introduced, and their 
finite dimensional representations were studied. The 
representation theory of these algebras was further 
studied in \cite{T1} from the perspectives of both 
spectral theory \cite{R} and Whittaker model \cite{K}. 
In \cite{JWS}, as generalizations of the algebras $U_{q}(f(K))$, 
another general class of algebras $U_{q}(f(K,H))$ was introduced 
and studied. Note that the Drinfeld quantum double of the 
positive part of the quantized enveloping algebra $U_{q}(sl_{2})$ 
studied in \cite{HZ} or equivalently the two-parameter 
quantum groups $U_{r,s}(sl_{2})$ studied in \cite{BW} is 
a special case of the algebra $U_{q}(f(K,H))$. The condition 
on the parameter Laurent polynomial $f(K,H)\in \C[K^{\pm 1}, H^{\pm
  1}]$ for the existence of a Hopf algebra structure on
$U_{q}(f(K,H))$ was determined, and finite dimensional 
irreducible representations were explicitly constructed 
as quotients of highest weight representations in \cite{JWS}. 
This class of algebras provides a family of quantum 
groups in the sense of Drinfeld \cite{D}. In particular,
$U_{q}(f_{m}(K,H))$ are quantum groups for $f_{m}(K,H) 
=\frac{K^{m}-H^{m}}{q-q^{-1}}, m\in \N$.

In this paper, we study the irreducible representations 
of these quantum groups $U_{q}(f_{m}(K,H))$. Though most of the 
results in this paper hold for the algebras $U_{q}(f(K,H))$ 
with general parameters, the calculations are more complicated. 

It is not surprising that these quantum groups share many 
similar properties with the two-parameter quantum groups 
$U_{r,s}(sl_{2})$. However, it may be useful to get more 
explicit information on the representation theory of these 
quantum groups. In the first part of this paper, we study the 
irreducible weight representations (which are not necessarily 
finite dimensional) of $U_{q}(f_{m}(K,H))$ from the viewpoint 
of spectral theory. Namely, we realize these quantum 
groups as Hyperbolic algebras, then apply the general results 
on Hyperbolic algebras established in \cite{R} to construct 
natural families of irreducible weight representations for 
$U_{q}(f_{m}(K,H))$. Such an approach yields the highest 
weight, the lowest weight and weight irreducible 
representations for $U_{q}(f_{m}(K,H))$. 

We denote by $f_{m}(K)$ the Laurent polynomial
$\frac{K^{m}-K^{-m}}{q-q^{-1}}$. Note that there is a close 
relationship between the representation theory of $U_{q}(f_{m}(K))$ 
and that of $U_{q}(f_{m}(K,H))$. We investigate this relationship 
following the idea in \cite{HZ}. As an application, we obtain some 
nice results on the category of all weight representations of 
$U_{q}(f_{m}(K,H))$. In particular, we show that the category 
of all weight representations of $U_{q}(f_{m}(K,H))$ is equivalent 
to the product of the category of weight representations of 
$U_{q}(f_{m}(K))$ with $\C^{\ast}$ as a tensor category. 
Combined with a result proved for $U_{q}(f_{m}(K))$ in \cite{JWZ}, 
we show that any finite dimensional weight representation of 
$U_{q}(f_{m}(K,H))$ is completely reducible. 

Finally, we study the Whittaker model for the center of these 
quantum groups. We prove that any Whittaker representation is irreducible 
if and only if it admits a central character. This criterion gives 
a complete classification of all irreducible Whittaker representations 
of $U_{q}(f_{m}(K,H))$.

The paper is organized as follows. In Section 1, we recall the
definitions of $U_{q}(f(K))$ and $U_{q}(f(K,H))$, and some basic 
facts about them from \cite{JWS,JWZ}. In Section 2, 
we recall some basic facts about spectral theory and Hyperbolic 
algebras from \cite{R}. Then we realize $U_{q}(f(K,H))$ as Hyperbolic 
algebras, and construct natural families of irreducible weight 
representations for $U_{q}(f_{m}(K,H))$. In Section 3, we study 
the relationship between $U_{q}(f_{m}(K))$ and $U_{q}(f_{m}(K,H))$ 
from the perspective of representation theory. In Section 4, we 
construct the Whittaker model for the center of $U_{q}(f_{m}(K,H))$, 
and study the Whittaker representations of $U(f_{m}(K,H))$. We obtain 
a classification of all irreducible Whittaker representations. 

\section{The algebras $U_{q}(f(K, H))$}
Let $\mathbb{C}$ be the field of complex numbers and $0 \neq q \in
\mathbb{C}$ such that $q^{2}\neq 1$. It is well-known that the 
quantized enveloping algebra $U_{q}(sl_{2})$ corresponding to 
the simple Lie algebra $sl_{2}$ is the associative $\C-$algebra 
generated by $K^{\pm 1}, E, F$ subject to the following relations:
\[
KE=q^{2}EK,\quad KF=q^{-2}FK,\quad KK^{-1}=K^{-1}K=1,
\]
\[
EF-FE=\frac{K-K^{-1}}{q-q^{-1}}.
\]

Note that $U_{q}(sl_{2})$ is a Hopf algebra with a Hopf algebra 
structure defined as follows:
\[
\Delta(E)=E\otimes 1+ K\otimes E,\quad \Delta(F)=F\otimes K^{-1}+1\otimes F;\]
\[
\epsilon(E)=0=\epsilon(F),\quad \epsilon(K)=1=\epsilon(K^{-1});
\]
\[
s(E)=-K^{-1}E,\quad s(F)=-FK,\quad s(K)=K^{-1}.
\]

As generalizations of $U_{q}(sl_{2})$, a class of algebras
$U_{q}(f(K))$ parameterized by Laurent polynomials $f(K)\in
\mathbb{C}[K,K^{-1}]$ was introduced in \cite{JWZ}. For the 
reader's convenience, we recall their definition here.
\begin{defn}
(See \cite{JWZ}) For any Laurent polynomial $f(K)\in \mathbb{C}[K,K^{-1}]$, 
$U_{q}(f(K))$ is defined to be the $\mathbb{C}-$algebra generated by 
$E,\,F,\,K^{\pm 1}$ subject to the following relations:
\[
KE=q^{2}EK,\quad KF=q^{-2}FK;
\]
\[
KK^{-1}=K^{-1}K=1;
\]
\[
EF-FE=f(K).
\]
\end{defn}

The ring theoretic properties and finite dimensional representations 
were first studied in detail in \cite{JWZ}. We state some of these 
results here without proof. First of all, for the Laurent polynomials 
$f(K)=a(K^{m}-K^{-m})$ where $a\in \C^{\ast}$ and $m\in \N$, the 
algebras $U_{q}(f(K))$ have a Hopf algebra structure. In particular, 
we have the following result from \cite{JWZ}:
\begin{prop}
(Prop 3.3 in \cite{JWZ}) Assume $f(K)$ is a non-zero Laurent 
polynomial in $\mathbb{C}[K,K^{-1}]$. Then the non-commutative 
algebra $U_{q}(f(K))$ is a Hopf algebra such that 
$K,K^{-1}$ are group-like elements, and $E, F$ are 
skew primitive elements if and only if $f(K)=a(K^{m}-K^{-m})$ 
with $m=t-s$ and the following conditions are satisfied:
\begin{gather*}
\Delta(K)=K\otimes K,\quad \Delta(K^{-1})=K^{-1}\otimes K^{-1};\\
\Delta(E)=E^{s}\otimes E + E\otimes K^{t},\quad \Delta(F)=K^{-t}\otimes F+F\otimes K^{-s};\\
\epsilon(K)=\epsilon(K^{-1})=1,\quad \epsilon(E)=\epsilon(F)=0;\\
S(K)=K^{-1},\quad S(K^{-1})=K;\\
S(E)=-K^{-s}EK^{-t},\quad S(F)=-K^{t}FK^{s}.
\end{gather*}
\end{prop}
\qed

For the case $f_{m}(K)=\frac{K^{m}-K^{-m}}{q-q^{-1}}$ for $m\in \mathbb{N}$
and $q$ is not a root of unity, the finite dimensional irreducible
representations were proved to be highest weight and constructed explicitly in
\cite{JWZ}. Furthermore, any finite dimensional representations are
completely reducible as stated in the following theorem from \cite{JWZ}.
\begin{thm}
(Thm 4.17 in \cite{JWZ}) With the above assumption for $f_{m}(K)$ and $q$, 
any finite dimensional representation $V$ of $U_{q}(f_{m}(K))$ is
completely reducible.
\end{thm}\qed

\begin{rem}
The representation theory of $U_{q}(f_{m}(K))$ was studied further from
the points of views of spectral theory and Whittaker model in \cite{T1}, 
where more families of interesting irreducible representations 
were constructed.
\end{rem}

As generalizations of the algebras $U_{q}(f(K))$, another general
class of algebras parameterized by Laurent polynomials $f(K,H)\in
\C[K^{\pm 1}, H^{\pm 1}]$ was introduced and studied in
\cite{JWS}. First, let us recall the definition of $U_{q}(f(K,H))$ here:
\begin{defn}
(See \cite{JWS}) Let $f(K, H)\in \mathbb{C}[K^{\pm 1},H^{\pm 1}]$ be
a Laurent polynomial, $U_{q}(f(K,H))$ is defined to be the
$\mathbb{C}-$algebra generated by $E,\,F,\,K^{\pm 1}, H^{\pm 1}$ subject to 
the following relations:
\begin{gather*}
KE=q^{2}EK,\quad KF=q^{-2}FK,\\
HE=q^{-2}EH,\quad HF=q^{2}FH,\\
KK^{-1}=K^{-1}K=1=HH^{-1}=H^{-1}H, \quad KH=HK,\\
EF-FE=f(K, H).
\end{gather*}
\end{defn}

It is easy to see that the Drinfeld quantum double of the 
positive part of $U_{q}(sl_{2})$ \cite{HZ} or the two-parameter 
quantum group $U_{r,s}(sl_{2})$ \cite{BW} is a special case of 
the algebra $U_{q}(f(K,H))$. The condition on the parameter $
f(K,H)$ for the existence of a Hopf algebra structure on 
$U_{q}(f(K,H))$ was determined, and finite dimensional irreducible 
representations were constructed explicitly as quotients of highest
weight representations in \cite{JWS}. In addition, a counter example 
was also constructed to show that not all finite dimensional representations
are completely reducible in \cite{JWS}. So it would be interesting to know what
kind of finite dimensional representations are completely
reducible. We will address this question in Section 3.

\section{Hyperbolic algebras and their representations}
In this section, we realize $U_{q}(f(K, H))$ as Hyperbolic algebras 
and apply the methods in spectral theory as developed in \cite{R} to 
construct irreducible weight representations of $U_{q}(f_{m}(K, H))$. For
the reader's convenience, we need to recall a little bit of background 
about spectral theory and Hyperbolic algebras from \cite{R}.

\subsection{Preliminaries on spectral theory}

Spectral theory of abelian categories was first started by Gabriel in
\cite{G}. He defined the injective spectrum of any noetherian Grothendieck 
category. This spectrum consists of isomorphism classes of
indecomposable injective objects. If $R$ is a commutative noetherian
ring, then the spectrum of the category of all $R-$modules is
isomorphic to the prime spectrum $Spec(R)$ of $R$. And one can 
reconstruct any noetherian commutative scheme $(X,O_{X})$ using the 
spectrum of the category of quasi-coherent sheaves of modules on $X$. 
The spectrum of any abelian category was later on defined by Rosenberg 
in \cite{R}. This spectrum works for any abelian category. Via this 
spectrum, one can reconstruct any quasi-separated and 
quasi-compact commutative scheme $(X,O_{X})$ via the spectrum of the 
category of quasi-coherent sheaves of modules on $X$. 

To proceed, we review some basic notions and facts about 
spectrum of any abelian category. First of all, we recall 
the definition of the spectrum of any abelian category, 
then we explain its applications in representation theory. 
We refer the reader to \cite{R} for more details.

Let $C_{X}$ be an abelian category and $M, N \in C_{X}$ 
be any two objects; We say that $M \succ N $ if and only 
if $N$ is a sub-quotient of the direct sum of finitely 
many copies of $M$. It is easy to verify that $\succ$  is a pre-order. 
We say $M \approx N$ if and only if $M\succ N$ and $N\succ M$. 
It is obvious that $\approx$ is an equivalence. Let $Spec(X)$ 
be the family of all nonzero objects $M\in C_{X}$ such that 
for any non-zero sub-object $N$ of $M$, $N\succ M$. 
\begin{defn}
(See \cite{R}) The spectrum of any abelian category is defined to be: 
\[
{\bf Spec(X)}=Spec(X)/\approx.
\]
\end{defn}

Though spectral theory is more important for the purpose of 
non-commutative algebraic geometry, it has nice applications 
to representation theory. The notion of the spectrum has a natural 
analogue of the Zariski topology. Under certain mild finiteness 
conditions, its closed points are in a one-to-one correspondence 
with the irreducible objects of the category. In particular, 
this is true for the category of representations of an algebra. 
To study irreducible representations, one can study the spectrum 
of the category of all representations, then single out closed 
points of the spectrum with respect to the associated topology. 

\subsection{The left spectrum of a ring}
If $C_X$ is the category $A-mod$ of left modules over a ring $A$, then 
it is sometimes convenient to express the points of ${\bf Spec}(X)$ in 
terms of left ideals of the ring $A$. In order to do so, the 
{\it left spectrum} $Spec_l(A)$ was defined in \cite{R}, which 
is by definition the set of all left ideals $p$ of $A$ such 
that $A/p$ is an object of $Spec(X)$. The relation 
$\succ$ on $A-mod$ induces a {\it specialization} relation 
among left ideals, in particular, the specialization relation 
on $Spec_l(A)$. Namely, $A/m\succ A/n$ iff there exists a 
finite subset $x$ of elements of $A$ such that such that the 
ideal $(n:x)=\{a\in A\ |\ ax\subset n\}$ is contained in $m$. 
Following \cite{R}, we denote this by $n\le m$. Note that 
the relation $\le$ is just the inclusion if $n$ is a two-sided 
ideal. In particular, it is the inclusion if the ring $A$ is 
commutative. The map which assigns to an element of $Spec_l(A)$ induces 
a bijection of the quotient $Spec_l(A)/\approx$ of $Spec_l(A)$ 
by the equivalence relation associated with $\le$ onto 
${\bf Spec}(X)$. From now on, we will not distinguish 
$Spec_l(A)/\approx$ from ${\bf Spec}(X)$ and will express 
results in terms of the left spectrum. 

\subsection{Hyperbolic algebra $R\{\xi,\theta\}$ and its spectrum}
Hyperbolic algebras are studied by Rosenberg in \cite{R} and  by
Bavula under the name of Generalized Weyl algebras in \cite{B}. 
Hyperbolic algebra structure is very convenient for the construction 
of points of the spectrum. As an application of spectral theory 
to representation theory, points of the left spectrum have been 
constructed for Hyperbolic algebras in \cite{R}. Many `small' 
algebras including the first Weyl algebra $A_{1}$, the enveloping 
algebra $U(sl_{2})$, and their quantized versions or deformations 
are Hyperbolic algebras. We will review some basic facts about 
Hyperbolic algebras and two important construction theorems from \cite{R}.

Let $\theta$ be an automorphism of a commutative algebra $R$; and let 
$\xi$ be an element of $R$. 
\begin{defn}
The Hyperbolic algebra $R\{\theta,\xi\}$ is defined to be the $R-$algebra 
generated by $x,y$ subject to the following relations:
\[
xy=\xi,\quad yx=\theta^{-1}(\xi)
\]
and 
\[
xa=\theta(a)x,\quad ya=\theta^{-1}(a)y
\]
for any $a\in R$. $R\{\theta,\xi\}$ is called a Hyperbolic algebra over $R$.
\end{defn}

Let $C_{X}=C_{R\{\theta,\xi\}}$ be the category of modules over 
$R\{\theta,\xi\}$. We denote by ${\bf Spec(X)}$ the spectrum of 
$C_{X}$. Points of the left spectrum of Hyperbolic algebras are 
studied in \cite{R}, and in particular we have the following 
construction theorems from \cite{R}.
\begin{thm}
(Thm 3.2.2.in \cite{R})
\begin{enumerate}
\item Let $P\in Spec(R)$, and assume that the orbit of $P$ under the 
action of the $\theta$ is infinite.
\begin{enumerate}
\item If $\theta^{-1}(\xi)\in P$, and $\xi \in P$, then the left ideal 
\[
P_{1,1}\colon=P+R\{\theta,\xi\}x+R\{\theta,\xi\}y
\]
is a two-sided ideal from $Spec_{l}(R\{\theta,\xi\})$.

\item If $\theta^{-1}(\xi)\in P$, $\theta^{i}(\xi)\notin P$ for $0\leq
  i\leq n-1$, 
and $\theta^{n}(\xi)\in P$, then the left ideal 
\[
P_{1,n+1}\colon=R\{\theta,\xi\}P+R\{\theta,\xi\}x+R\{\theta,\xi\}y^{n+1}
\]
belongs to $Spec_{l}(R\{\theta,\xi\})$.

\item If $\theta^{i}(\xi)\notin P$ for $i\geq 0$ and
  $\theta^{-1}(\xi)\in P$, then
\[
P_{1,\infty}\colon=R\{\theta,\xi\}P+R\{\theta,\xi\}x
\]
belongs to $Spec_{l}(R\{\theta,\xi\})$.

\item If $\xi \in P $ and $\theta^{-i}(\xi)\notin P$ for all $i\geq
  1$, 
then the left ideal 
\[
P_{\infty,1}\colon=R\{\theta,\xi\}P+R\{\theta,\xi\}y
\]
belongs to $Spec_{l}(R\{\theta,\xi\})$.
\end{enumerate}
\item If the ideal $P$ in (b),\, (c)\, or (d) is maximal, 
then the corresponding left ideal of $Spec_{l}(R\{\theta,\xi\})$ is maximal.

\item Every left ideal $Q \in Spec_{l}(R\{\theta,\xi\})$ such that
  $\theta^{\nu}(\xi)\in Q$ for 
a $\nu \in \Z$ is equivalent to one left ideal as defined above
  uniquely from a prime ideal $P \in Spec(R)$. The latter means that 
if $P$ and $P'$ are two prime ideals of $R$ and $(\alpha,\beta)$ and
  $(\nu,\mu)$ 
take values $(1,\infty),(\infty,1),(\infty,\infty)$ or $(1,n)$, 
then $P_{\alpha,\beta}$ is equivalent to $P'_{\nu,\mu}$ if and only 
if $\alpha=\nu,\beta=\mu$ and $P=P'$.
\end{enumerate}
\end{thm}
\qed

\begin{thm}
(Prop 3.2.3. in \cite{R}) 
\begin{enumerate}
\item Let $P\in Spec(R)$ be a prime ideal of $R$ 
such that $\theta^{i}(\xi)\notin P$ for $i\in \Z$ and $\theta^{i}(P)-P\neq \O$ 
for $i\neq 0$, then $P_{\infty,\infty}=R\{\xi,\theta\}P\in Spec_{l}(R\{\xi,\theta\})$.
\item Moreover, if ${\bf P}$ is a left ideal of $R\{\theta,\xi\}$ such 
that ${\bf P}\cap R=P$, then ${\bf P}=P_{\infty,\infty}$. In particular,
if $P$ is a maximal ideal, then $P_{\infty,\infty}$ is a maximal 
left ideal.
\item If a prime ideal $P'\subset R$ is such that
  $P_{\infty,\infty}=P'_{\infty,\infty}$, 
then $P'=\theta^{n}(P)$ for some integer $n$. Conversely,
$\theta^{n}(P)_{\infty,\infty}=P_{\infty,\infty}$ for any $n\in \Z$. 
\end{enumerate}
\end{thm}\qed

\subsection{Realize $U_{q}(f(K, H))$ as Hyperbolic algebras}
Let $R$ be the sub-algebra of $U_{q}(f(K,H))$ generated by 
$EF,\,K^{\pm 1}, \,H^{\pm 1}$, then $R$ is a commutative algebra. We define 
an algebra automorphism $\theta \colon R \longrightarrow R$ of $R$ by setting 
\begin{gather*}
\theta(EF)=EF+f(\theta(K),\theta(H)), \\
\theta(K^{\pm 1})=q^{\mp 2}K^{\pm 1},\\
\theta(H^{\pm 1})=q^{\pm 2}H^{\pm 1}.
\end{gather*}

It is easy to see that $\theta$ extends to an algebra automorphism of $R$.
Furthermore, we have the following lemma:
\begin{lem} The following identities hold:
\begin{gather*}
E(EF)=\theta(EF)E,\\
F(EF)=\theta^{-1}(EF)F,\\
EK=\theta(K)E,\\
FK=\theta^{-1}(K)F,\\
EH=\theta(H)E,\\
FH=\theta^{-1}(H)F.
\end{gather*}
\end{lem}
{\bf Proof:} We only verify the first one and the rest of them can be
checked similarly. 
\begin{eqnarray*}
E(EF)&=&E(FE+f(K,H))\\
&=&(EF)E+Ef(K,H)\\
&=&(EF)E+f(\theta(K),\theta(H))E\\
&=&(EF+f(\theta(K),\theta(H)))E\\
&=&\theta(EF)E.
\end{eqnarray*}
So we are done.\qed

From Lemma 2.1, we have the following result:
\begin{prop}
$U_{q}(f(K, H))=R\{\xi=EF,\theta\}$ is a Hyperbolic algebra
 with $R$ and $\theta$ defined as above.
\end{prop}
\qed

It easy to see that we have the following corollary:
\begin{cor}
(See also Prop. 2.5 in \cite{JWS}) $U_{q}(f(K,H))$ is noetherian 
domain of GK-dimension $4$.
\end{cor}
\qed
\subsection{Families of irreducible weight representations of
  $U_{q}(f_{m}(K, H))$} 
Now we can apply the above construction theorems to the case 
of $U_{q}(f_{m}(K, H))$, and construct families of irreducible weight 
representations of $U_{q}(f_{m}(K, H))$. 

Given $\alpha, \beta, \gamma \in \C$, we denote by 
\[
M_{\alpha,\beta,\gamma}=(\xi-\alpha, K-\beta, H-\gamma )\subset R
\]
the maximal ideal of $R$ generated by $\xi-\alpha, K-\beta,
H-\gamma$. We have the following lemma:
\begin{lem}
$\theta^{n}(M_{\alpha,\beta,\gamma})\neq M_{\alpha,\beta,\gamma}$ for 
any $n \geq 1$. In particular, $M_{\alpha,\beta,\gamma}$ has an
infinite orbit under the action of $\theta$.
\end{lem}

{\bf Proof:} We have 
\begin{eqnarray*}
\theta^{n}(K-\beta) & = & (q^{-2n}K-\beta)\\
& =& q^{-2n}(K-q^{2n}\beta).
\end{eqnarray*}
Since $q$ is not a root of unity, $q^{2n}\neq 1 $ for any $n\neq 0$.
So we have $\theta^{n}(M_{\alpha,\beta,\gamma})\neq M_{\alpha,\beta,\gamma}$ for any
$n\geq 1$.\qed

Now we construct all irreducible weight representations of $U_{q}(f_{m}(K,
H))$ with $f_{m}(K,H)=\frac{K^{m}-H^{m}}{q-q^{-1}}, m\in \N$. 

First of all, another lemma is in order:
\begin{lem} For $n\geq 0$, we have the following:
\begin{enumerate}
\item $
\theta^{n}(EF)=EF+\frac{1}{q-q^{-1}}(\frac{q^{-2m}(1-q^{-2nm})}{1-q^{-2m}}K^{m}-\frac{q^{2m}(1-q^{2nm})}{1-q^{2m}}H^{m}).$
\item $
\theta^{-n}(EF)=EF-\frac{1}{q-q^{-1}}(\frac{1-q^{2nm}}{1-q^{2m}}K^{m}-\frac{1-q^{-2nm}}{1-q^{-2m}}H^{m}).$
\end{enumerate}
\end{lem}

{\bf Proof:} For $n\geq 0$, we have 
\begin{eqnarray*}
\theta^{n}(EF)&=&EF+\frac{1}{q-q^{-1}}((q^{-2m}+\cdots+q^{-2nm})K^{m}\\
&-&(q^{2m}+\cdots+q^{2nm})H^{m})\\
&=& EF+\frac{1}{q-q^{-1}}(\frac{q^{-2m}(1-q^{-2nm})}{1-q^{-2m}}K^{m}-\frac{q^{2m}(1-q^{2nm})}{1-q^{2m}}H^{m}).
\end{eqnarray*}
The second statement can be verified similarly.\qed

\begin{thm} Let $P=M_{\alpha,\beta, \gamma}$, then we have the following:
\begin{enumerate}
\item If $\alpha= \frac{\beta^{m}-\gamma^{m}}{q-q^{-1}},\,
(\beta/\gamma)^{m}=q^{2mn}$ for some $n\geq 0$, then 
$\theta^{n}(\xi)\in M_{\alpha,\beta,\gamma}$ and $\theta^{-1}(\xi) \in
M_{\alpha,\beta,\gamma}$, thus $U_{q}(f_{m}(K,H))/P_{1,n+1}$ is a finite dimensional 
irreducible representation of $U_{q}(f_{m}(K,H))$.

\item If $\alpha= \frac{\beta^{m}-\gamma^{m}}{q-q^{-1}}$ and 
$(\beta/\gamma)^{m}\neq q^{2mn}$ for all $n\geq 0$, 
then $U_{q}(f_{m}(K,H))/P_{1,\infty}$ is an infinite 
dimensional irreducible representation of $U_{q}(f_{m}(K,H))$.

\item If $\alpha=0 $ and $0\neq \frac{1}{q-q^{-1}}(\frac{1-q^{2nm}}{1-q^{2m}}\beta^{m}-\frac{1-q^{-2nm}}{1-q^{-2m}}\gamma^{m})$ for any $n\geq 1$, 
then $U_{q}(f_{m}(K,H))/P_{\infty,1}$ is an infinite dimensional 
irreducible representation of $U_{q}(f_{m}(K,H))$. 
\end{enumerate}
\end{thm}
{\bf Proof:} Since
$\theta^{-1}(\xi)=\xi-\frac{K^{m}-H^{m}}{q-q^{-1}}$, 
thus $\theta^{-1}(\xi)\in M_{\alpha,\beta,\gamma}$ 
if and only if $\alpha=\frac{\beta^{m}-\gamma^{m}}{q-q^{-1}}$. 
Now by the proof of Lemma 2.3, we have 
\begin{eqnarray*} 
\theta^{n}(\xi)&=&\xi+\frac{1}{q-q^{-1}}((q^{-2m}+\cdots +q^{-2nm})K^{m}\\
&-&(q^{2m}+\cdots +q^{2nm})H^{m})\\
&=& \xi+\frac{1}{q-q^{-1}}(\frac{q^{-2m}(1-q^{-2nm})}{1-q^{-2m}}K^{m}-\frac{q^{2m}(1-q^{2nm})}{1-q^{2m}}H^{m}).
\end{eqnarray*}
Hence $\theta^{n}(\xi)\in M_{\alpha,\beta,\gamma}$ if and only if 
\begin{eqnarray*}
0&=&\alpha+\frac{1}{q-q^{-1}}((q^{-2m}+\cdots,+q^{-2nm})\beta^{m}-(q^{2m}+\cdots,+q^{2nm})\gamma^{m})\\
&=& \alpha+\frac{1}{q-q^{-1}}(\frac{q^{-2m}(1-q^{-2nm})}{1-q^{-2m}}\beta^{m}-\frac{q^{2m}(1-q^{2nm})}{1-q^{2m}}\gamma^{m}).
\end{eqnarray*}
Hence when $\alpha= \frac{\beta^{m}-\gamma^{m}}{q-q^{-1}},\,
(\beta/\gamma)^{m}=q^{2mn}$ for some $n\geq 0$, we have
\[ 
\theta^{n}(\xi)\in M_{\alpha,\beta,\gamma},\, \theta^{-1}(\xi) \in M_{\alpha,\beta,\gamma}.
\]
Thus by Theorem 2.1, $U_{q}(f_{m}(K,H))/P_{1,n+1}$ is a finite dimensional 
irreducible representation of $U_{q}(f_{m}(K,H))$.
So we have already proved the first statement, the rest of the statements
can be similarly verified.\qed

\begin{rem}
The representations we constructed in Theorem 2.3 exhaust all finite
dimensional irreducible weight representations, the highest weight irreducible
representations and the lowest weight irreducible representations of $U_{q}(f_{m}(K,H))$.
\end{rem}
\begin{rem}
Finite dimensional irreducible weight representations have been
constructed in \cite{JWS} as quotients of highest weight
representations. And a counter example has also been constructed in
\cite{JWS} to indicate that not all finite dimensional representations are
completely reducible.
\end{rem}
Apply the second construction theorem, we have the following theorem:
\begin{thm}
Let $P=M_{\alpha,\beta,\gamma}$. If  $\alpha \neq -\frac{1}{q-q^{-1}}(\frac{q^{-2m}(1-q^{-2nm})}{1-q^{-2m}}\beta^{m}-\frac{q^{2m}(1-q^{2nm})}{1-q^{2m}}\gamma^{m})$ for any $n\geq 0$ and
 $\alpha \neq \frac{1}{q-q^{-1}}(\frac{1-q^{2nm}}{1-q^{2m}}\beta^{m}-                                   
\frac{1-q^{-2nm}}{1-q^{-2m}}\gamma^{m})$ for any $n\geq 1$, then
 $U_{q}(f_{m}(K,H))/P_{\infty,\infty}$ 
is an infinite dimensional irreducible weight representation of $U_{q}(f_{m}(K,H))$.
\end{thm}
{\bf Proof:} The proof is very similar to that of Theorem 2.3, we will omit it here.\qed
\begin{cor}
The representations constructed in Theorem 2.3 and Theorem 2.4 exhaust all
irreducible weight representations of $U_{q}(f_{m}(K,H))$.
\end{cor}
{\bf Proof:} It follows directly from Theorems 2.1, 2.2, 2.3 and 2.4.\qed

\section{The relationship between $U_{q}(f_{m}(K))$ and $U_{q}(f_{m}(K,H))$}
Recall that we denote by $f_{m}(K,H)$ the polynomial
$\frac{K^{m}-H^{m}}{q-q^{-1}}$, and by $f_{m}(K)$ the Laurent 
polynomial $\frac{K^{m}-K^{-m}}{q-q^{-1}}$. We compare the quantum
groups $U_{q}(f_{m}(K))$ and $U_{q}(f_{m}(K,H))$. As a
result, we prove that any finite dimensional weight 
representation of $U_{q}(f_{m}(K, H))$ is completely reducible.

First of all, it is easy to see that we have the following lemma
generalizing the situation in \cite{HZ}:
\begin{lem}
The map which sends $E$ to $E$, $F$ to $F$, $K^{\pm 1}$ to $K^{\pm 1}$
, and $H^{\pm 1}$ to $K^{\mp 1}$ extends to a unique surjective Hopf
algebra homomorphism $\pi \colon U_{q}(f_{m}(K,H))\longrightarrow U_{q}(f_{m}(K))$.
\end{lem}
{\bf Proof:} Note that both $U_{q}(f_{m}(K))$ and $U_{q}(f_{m}(K,H))$
are Hopf algebras under the assumption on $f_{m}(K)$ and
$f_{m}(K,H)$. Since the kernel of $\pi$ is generated by $K-H^{-1}$, 
it is a Hopf ideal of $U_{q}(f_{m}(K,H))$. So we are done.\qed

Our goal in this section is to describe those representations $M$
of $U_{q}(f_{m}(K,H))$ such that $End_{U_{q}(f_{m}(K,H))}(M)=\C$. 
Since $KH$ is in the center and invertible, it acts on these 
representations by a non-zero scalar. As in \cite{HZ}, for 
each $z\in \C^{\ast}$, we define a $\C-$algebra homomorphism 
$\pi_{z}\colon U_{q}(f_{m}(K,H))\longrightarrow U_{q}(f_{m}(K))$ as 
follows:
\begin{gather*}
\pi_{z}(E)=z^{\frac{m}{2}}E,\quad \pi_{z}(F)=F;\\
\pi_{z}(K)=z^{\frac{1}{2}}K,\quad \pi_{z}(H)=z^{\frac{1}{2}}K^{-1}.
\end{gather*}

It is easy to see that $\pi_{z}$ is an algebra epimorphism with the
kernel of $\pi_{z}$ being a two-sided ideal generated by $KH-z$. But 
they may not necessarily be a Hopf algebra homomorphism unless $z=1$. 

Let $M$ be a representation of $U_{q}(f_{m}(K,H))$. As in \cite{HZ}, 
we have the following lemma:
\begin{lem}
Suppose that $End_{U_{q}(f_{m}(K,H))}(M)=\C$. Then there exists a unique 
$z\in \C^{\ast}$ such that $M$ is the pullback of a 
representation of $U_{q}(f_{m}(K))$ through a $\pi_{z}$ 
as defined above. In particular, any such irreducible 
representation of $U_{q}(f_{m}(K,H))$ is the pullback 
of an irreducible representation of $U_{q}(f_{m}(K))$ through 
the algebra homomorphism $\pi_{z}$ for some $z\in \C^{\ast}$.
\end{lem}\qed

We use the notation in \cite{HZ}. Let $M$ be a representation of 
$U_{q}(f_{m}(K))$, we denote by $M_{z}$ the representation of
$U_{q}(f_{m}(K,H))$ obtained as the pullback of $M$ via $\pi_{z}$. 
Let $\epsilon_{z}$ be the one dimensional representation of 
$U_{q}(f_{m}(K,H))$ which is defined by mapping the generators 
of $U_{q}(f_{m}(K,H))$ as follows:
\begin{gather*}
\epsilon_{z}(E)=\epsilon_{z}(F)=0;\\
\epsilon_{z}(K)=z^{\frac{1}{2}}, \quad \epsilon_{z}(H)=z^{\frac{1}{2}}.
\end{gather*}
Then we have the following similar lemma as in \cite{HZ}:
\begin{lem}
Let $0\neq z\in \C$, and $M$ be a representation of $U_{q}(f_{m}(K))$. Then
$M_{z}\cong \epsilon_{z}\otimes M_{1}\cong M_{1}\otimes
\epsilon_{z}$. In particular, if $0\neq z^{\prime}\in \C$ and $N$ is another
representation of $U_{q}(f_{m}(K))$, then we have 
\[
M_{z}\otimes N_{z'}\cong (M\otimes N)_{zz'}.
\]
\end{lem}
{\bf Proof:} The proof is straightforward. \qed

Let $M$ be a representation of $U_{q}(f_{m}(K,H))$. We say $M$ is a 
weight representation if $H$ and $K$ are acting on $M$ semisimply. 
Let $\frak C$ be the category of all weight representations of
$U_{q}(f_{m}(K))$ and $\tilde{\frak C}$ be the category of all weight 
representations of $U_{q}(f_{m}(K,H))$. Let $\C^{\ast}$ be the tensor 
category associated to the multiplicative group $\C^{\ast}$, then we 
have the following:
\begin{thm}
The category $\tilde{\frak C}$ is equivalent to the direct product 
of the categories $\frak C$ and $\C^{\ast}$ as a tensor category.
\end{thm}
{\bf Proof:} The proof is similar to the one in \cite{HZ}, we refer
the reader to \cite{HZ} for more details. \qed

\begin{cor}
Any finite dimensional weight representation of $U_{q}(f_{m}(K,H))$ is
completely reducible.
\end{cor}
{\bf Proof:} This follows from the above theorem and the fact that any
finite dimensional representation of $U_{q}(f_{m}(K))$ is completely
reducible (as is proved in \cite{JWZ}).
\qed
\begin{cor} 
The tensor product of any two finite dimensional weight
representations of $U_{q}(f_{m}(K,H))$ is completely reducible.
\end{cor}
\qed
\begin{rem}
After the first draft of this paper was written, we have been kindly 
informed by J. Hartwig that the complete reducibility of finite
dimensional weight representations is also proved in his preprint \cite{JH} 
in a more general setting of Ambiskew polynomial rings via a different 
approach.
\end{rem}

\begin{rem}
It might be interesting to study the decomposition of the 
product of two finite dimensional irreducible weight representations.
\end{rem}

\begin{rem}
When $m=1$, the above results are obtained in \cite{HZ} for the
Drinfeld double of the positive part of $U_{q}(sl_{2})$, and
equivalently for the two-parameter quantum groups $U_{r,s}(sl_{2})$ 
in \cite{BW}.
\end{rem}

\section{The Whittaker model for the center $Z(U_{q}(f_{m}(K,H)))$}
Let $\frak g$ be a finite dimensional complex semisimple Lie algebra 
and $U(\frak g)$ be its universal enveloping algebra. The Whittaker 
model for the center of $U(\frak g)$ was first studied by Kostant 
in \cite{K}. The Whittaker model for the center $Z(U(\frak g))$ is 
defined via a non-singular character of the nilpotent Lie 
subalgebra $\frak n^{+}$ of $\frak g$. Using the Whittaker 
model, Kostant studied the structure of Whittaker modules 
of $U(\frak g)$ and many important results about Whittaker 
modules were obtained in \cite{K}. Later on, Kostant's 
idea was further generalized by Lynch in \cite{L} and 
by Macdowell in \cite{M} to the case of singular characters 
of $\frak n^{+}$ and similar results were proved to hold in 
these cases. 

The obstacle of generalizing the Whittaker model to the quantized 
enveloping algebra $U_{q}(\frak g)$ with $\frak g$ of higher ranks 
is that there is no non-singular character existing for the positive 
part $(U_{q}(\frak g))^{> 0}$ of $U_{q}(\frak g)$ because of the
quantum Serre relations. In order to overcome this difficulty, it 
was Sevostyanov who first realized to use the topological version 
$U_{h}(\frak g)$ over $C[[h]]$ of quantum groups. Using a family 
of Coxeter realizations $U^{s_{\pi}}_{h}(\frak g)$ of the quantum 
group $U_{h}(\frak g)$ indexed by the Coxeter elements $s_{\pi}$, 
he was able to prove Kostant's results for $U_{h}(\frak g)$ in 
\cite{S1}. However, in the simplest case of ${\frak g}=sl_{2}$, 
the quantum Serre relations are trivial, thus a direct approach 
should still work and this possibility has been worked out 
recently in \cite{O}.

In addition, it is reasonable to expect that the Whittaker model 
exists for most of the deformations of $U_{q}(sl_{2})$. In this 
section, we show that there is such a Whittaker model for the 
center of $U_{q}(f_{m}(K,H))$, and will study the Whittaker modules 
for $U_{q}(f_{m}(K,H))$. We obtain analogous results as in \cite{K} 
and \cite{O}. For the reader's convenience, we present all the 
details here. Following the convention in \cite{K}, we will 
use the term of Whittaker modules instead of Whittaker 
representations.
\subsection{The center $Z(U_{q}(f_{m}(K,H)))$ of $U_{q}(f_{m}(K,H))$}
In this subsection, we give a description of the center of 
$U_{q}(f_{m}(K,H))$. The center $Z(U_{q}(f(K,H)))$ was also 
studied in \cite{JWS} as well. As mentioned at the very 
beginning, we always assume $f_{m}(K,H)=\frac{K^{m}-H^{m}}{q-q^{-1}}$ 
and $q$ is not a root of unity. 

We define a Casimir element of $U_{q}(f_{m}(H,K))$ by setting: 
\[ 
\Omega=FE+\frac{q^{2m}K^{m}+H^{m}}{(q^{2m}-1)(q-q^{-1})}.
\]
We have the following proposition:
\begin{prop}
\begin{eqnarray*}
\Omega &=& FE+\frac{q^{2m}K^{m}+H^{m}}{(q^{2m}-1)(q-q^{-1})}\\
  &=& EF+\frac{K^{m}+q^{2m}H^{m}}{(q^{2m}-1)(q-q^{-1})}.
\end{eqnarray*}
\end{prop}
{\bf Proof:} Since $EF=FE+\frac{K^{m}-H^{m}}{q-q^{-1}}$, we have
\begin{eqnarray*}
\Omega &=& FE+\frac{q^{2m}K^{m}+H^{m}}{(q^{2m}-1)(q-q^{-1})}\\
  &=& EF-\frac{K^{m}-H^{m}}{q-q^{-1}}+
\frac{q^{2m}K^{m}+H^{m}}{(q^{2m}-1)(q-q^{-1})}\\
&=& EF+\frac{K^{m}+q^{2m}H^{m}}{(q^{2m}-1)(q-q^{-1})}.
\end{eqnarray*}
So we are done.\qed

In addition, we have the following lemma:
\begin{lem}
$\Omega$ is in the center of $U_{q}(f_{m}(K,H))$.
\end{lem}
{\bf Proof:} It suffices to show that $\Omega E=E\Omega,\Omega
F=F\Omega,\Omega K=K\Omega,\Omega H =H\Omega$.
We will only verify that $\Omega E=E\Omega$ and the rest of them 
are similar.
\begin{eqnarray*}
\Omega E &=&(FE+\frac{q^{2m}K^{m}+H^{m}}{(q^{2m}-1)(q-q^{-1})})E\\
&=& (EF-\frac{K^{m}-H^{m}}{q-q^{-1}}+\frac{q^{2m}K^{m}+K^{-m}}{(q^{2m}-1)(q-q^{-1})})E\\
&=& E(FE)+\frac{K^{m}+q^{2m}H^{m}}{(q-q^{-1})(q^{2m}-1)}E\\
&=& E(FE+\frac{q^{2m}K^{m}+H^{m}}{(q^{2m}-1)(q-q^{-1})})\\
&=& E\Omega.
\end{eqnarray*}
So we are done with the proof.\qed

In particular, we have the following description of the center 
$Z(U_{q}(f_{m}(K,H)))$ of $U_{q}(f_{m}(K,H))$.
\begin{prop}
(See also \cite{JWS}) $Z(U_{q}(f_{m}(K,H)))$ is the subalgebra of $U_{q}(f_{m}(K,H))$ generated by $\Omega, (KH)^{\pm 1}$. In particular, $Z(U_{q}(f_{m}(K,H)))$ is 
isomorphic to the localization $(\C[\Omega, KH])_{(KH)}$ of the polynomial
ring in two variables $\Omega, KH$.
\end{prop}
{\bf Proof:} 
By Lemma 3.1., we have $\Omega,(KH)^{\pm 1} \in Z(U_{q}(f_{m}(K,H)))$. Thus 
the subalgebra $\C[\Omega,(KH)^{\pm 1}]$ generated by $\Omega,
(KH)^{\pm 1}$ is contained 
in $Z(U_{q}(f_{m}(K,H)))$. So it suffices to prove 
the other inclusion $Z(U_{q}(f_{m}(K,H)))\subseteq
\mathbb{C}[\Omega,(KH)^{\pm 1}]$. 
Note that $U_{q}(f_{m}(K,H))=\bigoplus_{n\in \Z_{\geq 0}}U_{q}(f_{m}(K,H))_{n}$ 
where $U_{q}(f_{m}(K,H))_{n}$ is the $\C-$span of elements
$\{u\in U_{q}(f_{m}(K,H))\mid Ku=q^{2n}uK,\ Hu=q^{-2n}uH\}$. 
Suppose $x\in Z(U_{q}(f_{m}(K,H)))$, 
then $xK=Kx,xH=Hx$. Thus $x\in U_{q}(f_{m}(K,H))_{0}$, 
which is generated by $EF,K^{\pm 1},H^{\pm 1}$. By the 
definition of $\Omega$, we know that $U_{q}(f_{m}(K,H))_{0}$ 
is also generated by $\Omega,K^{\pm 1}, H^{\pm 1}$. 
Hence $x=\sum f_{i,j}(\Omega)K^{i}H^{j}$ where $f_{i,j}(\Omega)$ 
are polynomials in $\Omega$. Therefore 
\[
xE=\sum f_{i,j}(\Omega)K^{i}H^{j}E=\sum f_{i,j}(\Omega)q^{2i-2j}EK^{i}H^{j}=Ex,
\]
which forces $i=j$. So $x\in \mathbb{C}[\Omega, (KH)^{\pm 1}]$ as desired. 
So we have proved that $Z(U_{q}(f_{m}(K,H)))=\mathbb{C}[\Omega,(KH)^{\pm 1}]$.
\qed

\subsection{The Whittaker model for $Z(U_{q}(f_{m}(K,H)))$}
Now we construct the Whittaker model for $Z(U_{q}(f_{m}(K,H)))$ following 
the lines in \cite{K} and \cite{O}. In the rest of this subsection, we will
denote the parameter Laurent polynomial $\frac{K^{m}-H^{m}}{q-q^{-1}}$
by $f(K,H)$ instead of $f_{m}(K,H)$.  

First, we fix some notations. We denote by $U_{q}(E)$ the subalgebra
of $U_{q}(f(K,H))$ generated by $E$, by $U_{q}(F,K^{\pm 1},H^{\pm 1})$ 
the subalgebra of $U_{q}(f(K,H))$ generated by $F,K^{\pm 1},H^{\pm
  1}$. A non-singular character of the algebra $U_{q}(E)$ can be
defined as follows:
\begin{defn}
An algebra homomorphism $\eta\colon U_{q}(E)\longrightarrow \mathbb{C}$ 
is called a non-singular character of $U_{q}(E)$ if $\eta(E)\neq 0$.
\end{defn}

From now on, we will fix such a non-singular character of $U_{q}(E)$
and denote it by $\eta$. Following \cite{K}, we define the concepts of 
a Whittaker vector and a Whittaker module corresponding to the fixed 
non-singular character $\eta$.
\begin{defn}
Let $V$ be a $U_{q}(f(K,H))-$module, a vector $0\neq v\in V$ is called 
a Whittaker vector of type $\eta$ if $E$ acts on $v$ through the 
non-singular character $\eta$, i.e., $Ev=\eta(E)v$. If
$V=U_{q}(f(K,H))v$, then we call $V$ a Whittaker module 
of type $\eta$ and $v$ is called a cyclic Whittaker 
vector of type $\eta$. 
\end{defn}

The following decomposition of $U_{q}(f(K,H))$ is obvious:
\begin{prop}
$U_{q}(f(K,H))$ is isomorphic to $U_{q}(F,K^{\pm 1},H^{\pm
  1})\otimes_{\mathbb{C}} U_{q}(E)$ as a vector space and 
$U_{q}(f(K,H))$ is a free module over the subalgebra $U_{q}(E)$.
\end{prop}
\qed

Let us denote the kernel of $\eta \colon U_{q}(E)\longrightarrow \mathbb{C}$ 
by $U_{q,\eta}(E)$, and we have the following decompositions 
of $U_{q}(E)$ and $U_{q}(f(K,H))$.
\begin{prop}
We have $U_{q}(E)=\mathbb{C} \oplus U_{q,\eta}(E)$. In addition, 
\[
U_{q}(f(K,H))
\cong U_{q}(F,K^{\pm 1},H^{\pm 1})\oplus U_{q}(f(K,H))U_{q,\eta}(E).
\]
\end{prop}
{\bf Proof:} It is obvious that $U_{q}(E)=\mathbb{C} \oplus U_{q,\eta}(E)$. 
And we have 
\[
U_{q}(f(K,H))=U_{q}(F,K^{\pm 1},H^{\pm 1})\otimes(\C \oplus U_{q,\eta}(E)),
\]

thus 
\[
U_{q}(f(K,H))\cong U_{q}(F,K^{\pm 1},H^{\pm 1})\oplus U_{q}(f(K,H))U_{q,\eta}(E).
\]
So we are done.
\qed

Now we define a projection:
\[
\pi\colon U_{q}(f(K,H))\longrightarrow U_{q}(F,K^{\pm 1},H^{\pm 1})
\]
 from $U_{q}(f(K,H))$ onto $U_{q}(F,K^{\pm 1},H^{\pm 1})$ by taking 
the $U_{q}(F,K^{\pm 1},H^{\pm 1})-$component of any $u\in U_{q}(f(K,H))$. 
We denote the image $\pi(u)$ of $u\in U_{q}(f(K,H))$ 
by $u^{\eta}$ for short.
\begin{lem}
If $v \in Z(U_{q}(f(K,H)))$ and $u\in U_{q}(f(K,H))$, then 
we have $u^{\eta}v^{\eta}=(uv)^{\eta}$.
\end{lem}
{\bf Proof:} Let $v \in Z(U_{q}(f(K,H))),\, u\in U_{q}(f(K,H))$, then we have
\begin{eqnarray*}
 uv-u^{\eta}v^{\eta}&=&(u-u^{\eta})v+u^{\eta}(v-v^{\eta})\\
&=& v(u-u^{\eta})+u^{\eta}(v-v^{\eta}),
\end{eqnarray*}
which is in $U_{q}(f(K,H))U_{q,\eta}(E)$. 
Hence $(uv)^{\eta}=u^{\eta}v^{\eta}$.
\qed

By the definition of $\Omega$, we have the following description of 
$\pi(\Omega)$:
\begin{lem}
\[
\pi(\Omega)=\eta(E)F+\frac{q^{2m}K^{m}+H^{m}}{(q^{2m}-1)(q-q^{-1})}.
\]
\end{lem}\qed

\begin{prop}
The map 
\[
\pi\colon Z(U_{q}(f(K,H))\longrightarrow U_{q}(F,K^{\pm 1},H^{\pm 1})
\]
is an algebra isomorphism of $Z(U_{q}(f(K,H)))$ onto 
its image $W(F,K^{\pm 1},H^{\pm 1})$ in $U_{q}(F,K^{\pm 1},H^{\pm 1})$.
\end{prop}
{\bf Proof:} It follows from Lemma 4.2. that $\pi$ is a homomorphism 
of algebras. By Lemma 4.3, we have
\[
\pi(\Omega)=\eta(E)F+\frac{q^{2m}K^{m}+H^{m}}{(q^{2m}-1)(q-q^{-1})}
\]
with $\eta(E)\neq 0$. Note that $\pi(KH)=KH$. We show that $\pi$ is
injective. Suppose that $\pi(z)=0$ for some element $0\neq z\in 
Z(U_{q}(f(K,H)))$. Since $Z(U_{q}(f(K,H))=\C[\Omega,(KH)^{\pm 1}]$, 
we can write $z=\sum_{i=0}^{k}z_{i}(KH) \Omega^{i}$ where $z_{i}(KH)$ 
are non-zero Laurent polynomials in $\C[(KH)^{\pm 1}]$. Since 
$\pi(z)=0$, then $u_{k}(KH)(\eta(E))^{k}F^{k}=0$, which is a 
contradiction. So $\pi$ is an injection. Thus $\pi$ is an 
algebra isomorphism from $Z(U_{q}(f(K,H)))$ onto its 
image $W(F,K^{\pm 1},H^{\pm 1})$ in $U_{q}(F,K^{\pm 1},H^{\pm 1})$.
\qed

\begin{lem}
If $u^{\eta}=u$, then we have 
\[
u^{\eta}v^{\eta}=(uv)^{\eta}
\]
 for any $v\in U_{q}(f(K,H))$.
\end{lem}
{\bf Proof:} We have 
\begin{eqnarray*}
uv-u^{\eta}v^{\eta}&=&(u-u^{\eta})v+u^{\eta}(v-v^{\eta})\\
&=&u^{\eta}(v-v^{\eta}),
\end{eqnarray*} 
which is in $U_{q}(f(K,H))U_{q,\eta}(E)$. So we have 
\[
u^{\eta}v^{\eta}=(uv)^{\eta}
\]
for any $v\in U_{q}(f(K,H))$.
\qed

Let $\tilde{A}$ be the subspace of $U_{q}(f(K,H))$ 
spanned by $K^{\pm i}$ where $i\in \Z_{\geq 0}$. Then $\tilde{A}$ is a graded 
vector space with 
\[
\tilde{A}_{[n]}=\C K^{n}\oplus \C K^{-n}
\]
for $n\geq 1$, and 
\[
\tilde{A}_{[0]}=\C,
\] 
and 
\[
\tilde{A}_{[n]}=0
\]
for $n\leq -1$. 

We define a filtration of $U_{q}(F,K^{\pm 1}, H^{\pm 1})$ as follows:
\[
U_{q}(F,K^{\pm 1}, H^{\pm 1})_{[n]}=\bigoplus_{im+\mid j- k\mid
  \leq nm} U_{q}(F,K^{\pm 1},H^{\pm 1})_{i, j, k}
\]
with $U_{q}(F,K^{\pm 1},H^{\pm 1})_{i,j,k}$ being the vector space 
spanned by $F^{i}K^{j}H^{k}$.

We denote by 
\[
W(F,K^{\pm 1},H^{\pm 1})_{[p]}=\C[(KH)^{\pm 1}]-span\{1,\Omega^{\eta},\cdots,(\Omega^{\eta})^{p}\}
\]
for $q\geq 0$. It is easy to see that
\[
W(F,K^{\pm 1},H^{\pm 1})_{[p]}\subset W(F,K^{\pm 1},H^{\pm
  1})_{[p+1]},
\] 
and
\[ 
W(F,K^{\pm 1},H^{\pm 1})=\sum_{p\geq 0}W(F,K^{\pm 1},H^{\pm
1})_{[p]}.
\]
Note that $W(F,K^{\pm 1},H^{\pm 1})_{[p]}$ give a filtration of 
$W(F,K^{\pm 1},H^{\pm 1})$ which is compatible with the filtration 
of $U_{q}(F,K^{\pm 1},H^{\pm 1})$. In particular, we have
\[
W(F,K^{\pm 1},H^{\pm 1})_{[p]}=W(F,K^{\pm 1},H^{\pm 1})\cap U_{q}(F,K^{\pm 1},H^{\pm 1})_{[p]}
\]
for $p\geq 0$ via direct computations.

Now, we have a decomposition of $U_{q}(F,K^{\pm 1},H^{\pm 1})$ as follows:
\begin{thm}
$U_{q}(F,K^{\pm 1},H^{\pm 1})$ is free (as a right module) over 
$W(F,K^{\pm 1},H^{\pm 1})$. And the multiplication induces an 
isomorphism 
\[
\Phi\colon \tilde{A}\otimes W(F,K^{\pm 1},H^{\pm 1})\longrightarrow
U_{q}(F,K^{\pm 1}, H^{\pm 1})
\]
as right $W(F,K^{\pm 1},H^{\pm 1})-$modules. In particular, we have 
the following:
\[
\bigoplus_{p+lm=nm}\tilde{A}_{[p]}\otimes W(F,K^{\pm 1},H^{\pm 1})_{[l]}\cong
U_{q}(F,K^{\pm 1},H^{\pm 1})_{[n]}.
\]
\end{thm}
{\bf Proof:} Note that the map $\tilde{A}\times W(F,K^{\pm 1},H^{\pm
  1}) \longrightarrow U_{q}(F,K^{\pm 1},H^{\pm 1})$ is bilinear. So by 
the universal property of the tensor product, there is a map 
from $\tilde{A}\otimes W(F,K^{\pm 1},H^{\pm 1})$ into 
$U_{q}(F,K^{\pm 1},H^{\pm 1})$ defined by the multiplication. 
It is easy to check this map is a homomorphism of right 
$W(F,K^{\pm 1},H^{\pm 1})-$modules and is surjective 
as well. 

Now, it remains to show that the map is injective. Let 
$0\neq u\in \tilde{A}\otimes W(F,K^{\pm 1},H^{\pm 1})$ with $\Phi(u)=0$. 
We can write
\[
u=\sum_{i=0}^{N}a_{i}(K)\otimes b_{i}(KH)(\eta(E)F+\frac{q^{2m}K^{m}+H^{m}}{(q^{2m}-1)(q-q^{-1})})^{i}
\]
where $u_{i}(K)$ are nonzero Laurent polynomials 
in $\C[K^{\pm 1}]$ and $b_{i}(KH)$ are non-zero 
Laurent polynomials in $\C[(KH)^{\pm 1}]$. Since $\Phi(u)=0$, 
then by direct computations, we have 
\[
a_{N}(K)b_{N}(KH)(\eta(E))^{N}F^{N}=0.
\]
Thus $a_{N}(K)b_{N}(KH)(\eta(E))^{N}=0$, which is a contradiction. 
So we have proved that $\Phi$ is indeed an isomorphism of vector spaces. 

In addition, by counting the degrees of both sides, we also have 
\[
\bigoplus_{p+lm=nm}\tilde{A}_{[p]}\otimes W(F,K^{\pm 1},H^{\pm 1})_{[l]}\cong
U_{q}(F,K^{\pm 1},H^{\pm 1})_{[n]}. 
\]
Thus the theorem is proved.
\qed

Let $Y_{\eta}$ be the left $U_{q}(f(K,H))-$module 
defined by 
\[
Y_{\eta}=U_{q}(f(K,H))\otimes_{U_{q}(E)} \C_{\eta},
\]
where $\C_{\eta}$ is the one dimensional 
$U_{q}(E)-$module defined by the character 
$\eta$. It is easy to see that 
\[
Y_{\eta}\cong U_{q}(f(K,H))/U_{q}(f(K,H))U_{q,\eta}(E)
\]
is a Whittaker module with a cyclic vector denoted by $1_{\eta}$. 
Now we have a quotient map from $U_{q}(f(K,H))$ to $Y_{\eta}$ 
as follows:
\[
U_{q}(f(K,H))\longrightarrow Y_{\eta}\quad \text{is defined by} 
\quad u\mapsto u1_{\eta},
\]
for any $u\in U_{q}(f(K,H))$.

If $u\in U_{q}(f(K,H))$, then there is a $u^{\eta}$ which is the
unique element in $U_{q}(F,K^{\pm 1},H^{\pm 1})$ such that
$u 1_{\eta}=u^{\eta}1_{\eta}$. As in \cite{K}, we define the $\eta
-$reduced action of $U_{q}(E)$ on $U_{q}(F,K^{\pm 1},H^{\pm 1})$ 
as follows:
\[
x\bullet v=(xv)^{\eta}-\eta(x)v,
\]
where $x\in U_{q}(E)$ and $v\in U_{q}(F,K^{\pm 1},H^{\pm 1})$.
\qed
\begin{lem}
Let $u\in U_{q}(f(K,H))$ and $x\in U_{q}(E)$, 
we have
\[
x\bullet u^{\eta}=[x,u]^{\eta}.
\]
\end{lem}
{\bf Proof:} $[x,u]1_{\eta}=(xu-ux)1_{\eta}=(xu-\eta(x)u)1_{\eta}$. Hence
\[[x,u]^{\eta}=(xu)^{\eta}-\eta(x)u^{\eta}=(xu^{\eta})^{\eta}-\eta(x)u^{\eta}=x\bullet u^{\eta}.\]\qed

\begin{lem}
Let $x\in U_{q}(E)$, $u\in U_{q}(F,K^{\pm 1},H^{\pm 1})$, and 
$v\in W(E,K^{\pm 1},H^{\pm 1})$, then we have 
\[
x \bullet (uv)=(x\bullet u)v.
\]
\end{lem}
{\bf Proof:} Let $v=w^{\eta}$ for some $w\in Z(U_{q}(f(K,H))$, 
then $uv=uw^{\eta}=u^{\eta}w^{\eta}=(uw)^{\eta}$. 
Thus
\begin{eqnarray*}
x\bullet(uv)&= &x\bullet(uw)^{\eta}=[x,uw]^{\eta}\\
&=&([x,u]w)^{\eta}=[x,u]^{\eta}w^{\eta}\\
&=&(x\bullet u^{\eta})v\\
&=&(x\bullet u)v.
\end{eqnarray*}
So we are done.\qed

Let $V$ be a $U_{q}(f(K,H))-$module and let $U_{q,V}(f(K,H))$ be 
the annihilator of $V$ in $U_{q}(f(K,H))$. Then $U_{q,V}(f(K,H))$ 
defines a central ideal $Z_{V}\subset Z(U_{q}(f(K,H)))$ by setting 
$Z_{V}=U_{q,V}(f(K,H))\cap Z(U_{q}(f(K,H)))$. Suppose that $V$ is 
a Whittaker module with a cyclic Whittaker  vector $w$, we denote 
by $U_{q,w}(f(K,H))$ the annihilator of $w$ in $U_{q}(f(K,H))$. It 
is obvious that 
\[
U_{q}(f(K,H))U_{q,\eta}(E)+U_{q}(f(K,H))Z_{V}\subset U_{q,w}(f(K,H)).
\]

In the next theorem, we show that the reverse inclusion holds.
First of all, we need an auxiliary Lemma:
\begin{lem}
Let $X=\{v\in U_{q}(F,K^{\pm 1},H^{\pm 1})\mid (x \bullet v)w=0,x\in U_{q}(E)\}$.
Then
\[
X=\tilde{A}\otimes W_{V}(F,K^{\pm 1},H^{\pm 1})+W(F,K^{\pm 1},H^{\pm 1}),
\]
where $W_{V}(F,K^{\pm 1},H^{\pm 1})=(Z_{V})^{\eta}$. In fact, 
$U_{q,V}(F,K^{\pm 1},H^{\pm 1})\subset X$ and 
\[
U_{q,w}(F,K^{\pm 1},H^{\pm 1})=\tilde{A}\otimes W_{w}(F,K^{\pm 1},H^{\pm 1}),
\]
where $U_{q,w}(F,K^{\pm 1},H^{\pm 1})=U_{q,w}(f(K,H))\cap
U_{q}(F,K^{\pm 1},H^{\pm 1})$.
\end{lem}
{\bf Proof:} Let us denote by $Y=\tilde{A}\otimes W_{V}(F,K^{\pm
  1},H^{\pm 1})+W(F,K^{\pm 1},H^{\pm 1})$ where $W(F,K^{\pm 1},H^{\pm 1})=(Z(U_{q}(f(K,H))))^{\eta}$. 
Thus we need to verify $X=Y$. Let $v\in W(F,K^{\pm 1},H^{\pm 1})$, 
then $v=u^{\eta}$ for some $u\in Z(U_{q}(f(K,H)))$. 
So we have 
\begin{eqnarray*}
x\bullet v &=&x\bullet u^{\eta}\\
&=&[x,u]^{\eta}\\
&=&(xu)^{\eta}-\eta(x)u^{\eta}\\
&=&x^{\eta}u^{\eta}-\eta(x)u^{\eta}\\
&=&0.
\end{eqnarray*}
So we have $W(F,K^{\pm 1},H^{\pm 1})\subset X$.
Let $u\in Z_{V}$ and $v\in U_{q}(F,K^{\pm 1},H^{\pm 1})$. 
Then for any $x\in U_{q}(E)$, we have 
\[
x\bullet(vu^{\eta})=(x\bullet v)u^{\eta}.
\]
Since $u\in Z_{V}$, then $u^{\eta}\in U_{q,w}(f(K,H))$. 
Thus we have $vu^{\eta}\in X$, hence 
\[
\tilde{A}\otimes W_{V}(F,K^{\pm 1},H^{\pm 1})\subset X,
\]
which proves $Y\subset X$. Note that $\tilde{A_{[i]}}$ 
is the subspace of $\C[K^{\pm 1}]$ spanned by 
$K^{\pm i}$, and let $\overline{W_{V}(F,K^{\pm 1},H^{\pm 1})}$ be 
the complement of 
$W_{V}(F,K^{\pm 1},H^{\pm 1})$ in $W(F,K^{\pm 1},H^{\pm 1})$. 

Let us set 
\[
M_{i}=\tilde{A_{[i]}}\otimes \overline{W_{V}(F,K^{\pm 1},H^{\pm 1})},
\]
thus we have the following:
\[
U_{q}(F,K^{\pm 1},H^{\pm 1})=M\oplus Y,
\]
where $M=\sum_{i\geq 1}M_{i}$.
We show that $M\cap X\neq 0$. Let $M_{[k]}=\sum_{1\leq i\leq k}M_{i}$, 
then $M_{[k]}$ are a filtration of $M$. Suppose $n$ is the smallest
integer such that $X\cap M_{[n]}\neq 0$ and $0\neq y\in X\cap
M_{[n]}$. Then we have $y=\sum_{1\leq i\leq n}y_{i}$ where $y_{i}\in
\tilde{A_{i}}\otimes \overline{W_{V}(F, K^{\pm 1},H^{\pm
    1})}$. Suppose we have chosen $y$ in such a way that $y$ has the
fewest terms. By similar computations as in \cite{O}, 
we have $0\neq y-\frac{1}{\eta(E)(q^{-2n}-1)}E\bullet y \in X\cap M_{[n]}$
with fewer terms than $y$. This is a contradiction. So we have $X\cap M=0$.

Now we prove that $U_{q,w}(F,K^{\pm 1},H^{\pm 1})\subset X$. 
Let $u\in U_{q,w}(F, K^{\pm 1}, H^{\pm 1})$ and $x\in U_{q}(E)$, 
then we have $xuw=0$ and $uxw=\eta(x)uw=0$. 
Thus $[x,u]\in U_{q,w}(f(K,H))$, hence $[x,u]^{\eta}\in U_{q,w}(F,K^{\pm 1},H^{\pm 1})$. 
Since $u\in U_{q,w}(F,K^{\pm 1},H^{\pm 1})\subset U_{q,w}(E,F,K^{\pm
  1},H^{\pm 1})$, 
then $x\bullet u=[x,u]^{\eta}$. 
Thus $x\bullet u \in U_{q,w}(F,K^{\pm 1},H^{\pm 1})$. 
So $u\in X$ by the definition of $X$. 
Now we are going to prove the following:
\[
W(F,K^{\pm 1},H^{\pm 1})\cap U_{q,w}(F,K^{\pm 1},H^{\pm
  1})=W_{V}(F,K^{\pm 1},H^{\pm 1}).
\] 
In fact, $W_{V}(F,K^{\pm 1},H^{\pm 1})=(Z_{V}^{\eta})$ and 
$W_{V}(F,K^{\pm 1},H^{\pm 1})\subset U_{q,w}(F,K^{\pm 1},H^{\pm 1})$. 
So if $v\in W_{w}(F,K^{\pm 1},H^{\pm 1})\cap U_{q,w}(F,K^{\pm 1},H^{\pm 1})$, 
then we can uniquely write $v=u^{\eta}$ for $u\in Z(U_{q}(f(K,H)))$. 
Then $vw=0$ implies $uw=0$ and hence 
$u\in Z(U_{q}(f(K,H)))\cap U_{q,w}(F, K^{\pm 1},H^{\pm 1})$. 
Since $V$ is generated cyclically by $w$, we 
have proved the above statement.
Obviously, we have 
$U_{q}(f(K,H))Z_{V}\subset U_{q,w}(f(K,H))$. 
Thus we have 
$\tilde{A}\otimes W_{V}(F,K^{\pm 1},H^{\pm 1})\subset U_{q,w}(F,K^{\pm
  1},H^{\pm 1})$. Therefore, we have 

\[U_{q,w}(F,K^{\pm 1},H^{\pm 1})=\tilde{A}\otimes
W_{V}(F,K^{\pm 1},H^{\pm 1}).\]

So we have finished the proof.\qed
\begin{thm}
Let $V$ be a Whittaker module admitting a cyclic 
Whittaker vector $w$, then we have 
\[
U_{q,w}(f(K,H))=U_{q}(f(K,H))Z_{V}+U_{q}(f(K,H))U_{q,\eta}(E).
\]
\end{thm}
{\bf Proof:} It is obvious that 
\[
U_{q}(f(K,H))Z_{V}+U_{q}(f(K,H))U_{\eta}(E)
\subset U_{q,w}(f(K,H)).
\]
Let $u\in U_{q,w}(f(K,H))$, we show that 
$u\in U_{q}(f(K,H))Z_{V}+U_{q}(f(K,H))U_{q,\eta}(E)$. 
Let $v=u^{\eta}$, then it suffices to show 
that $v\in \tilde{A}\otimes W_{V}(F,K^{\pm 1},H^{\pm 1})$. 
But $v\in U_{q,w}(F,K^{\pm 1},H^{\pm 1})=\tilde{A}\otimes
W_{V}(F,K^{\pm 1},H^{\pm 1})$.
So we have proved the theorem.\qed
\begin{thm}
Let $V$ be any Whittaker module for $U_{q}(f(K,H))$, 
then the correspondence
\[
V \mapsto Z_{V}
\]
sets up a bijection between the set of all 
equivalence classes of Whittaker modules and 
the set of all ideals of $Z(U_{q}(f(K,H)))$.
\end{thm}
{\bf Proof:} Let $V_{i}, i=1,2$ be two Whittaker modules. 
If $Z_{V_{1}}=Z_{V_{2}}$, then clearly $V_{1}$ is 
equivalent to $V_{2}$ by the above Theorem. 
Now let $Z_{\ast}$ be an ideal of $Z(U_{q}(f(K,H)))$ 
and let $L=U_{q}(f(K,H))Z_{\ast}+U_{q}(f(K,H))U_{\eta}(E)$. 
Then $V=U_{q}(f(K,H))/L$ is a Whittaker module 
with a cyclic Whittaker vector $w=\bar{1}$. 
Obviously we have $U_{q,w}(f(K,H))=L$. So 
$L=U_{q,w}(f(K,H))=U_{q}(f(K,H))Z_{V}+U_{q}(f(K,H))U_{q,\eta}(E)$. 
This implies that 
\[
\pi(Z_{\ast})=\pi(L)=\pi(Z_{V}).
\]
Since $\pi$ is injective on $Z(U_{q}(f(K,H)))$, thus $Z_{V}=Z_{\ast}$. 
Thus we finished the proof.\qed

\begin{thm}
Let $V$ be an $U_{q}(f(K,H))-$module. 
Then $V$ is a Whittaker module if and only if 
\[
V\cong U_{q}(f(K,H))\otimes_{Z(U_{q}(f(K,H)))
\otimes U_{q}(E)}(Z(U_{q}(f(K,H)))/Z_{\ast})_{\eta}.
\]
In particular, in such a case the ideal $Z_{\ast}$ 
is uniquely determined to be $Z_{V}$.
\end{thm}
{\bf Proof:} If $1_{\ast}$ is the image 
of $1$ in $Z(U_{q}(f(K,H)))/Z_{\ast}$, 
then 
\[
Ann_{Z(U_{q}(f(K,H)))\otimes U_{q}(F)}(1_{\ast})
=U_{q}(E)Z_{\ast}+Z(U_{q}(f(K,H)))U_{q,\eta}(E)
\]
Thus the annihilator of $w=1\otimes 1_{\ast}$ 
is 
\[
U_{q,w}(f(K,H))=U_{q}(f(K,H))Z_{\ast}+U_{q}(f(K,H))U_{q,\eta}(E)
\] 
Then the result follows from the last theorem. \qed

\begin{thm}
Let $V$ be an $U_{q}(f(K,H))-$module with a cyclic 
Whittaker vector $w\in V$. Then any $v\in V$ is a 
Whittaker vector if and only if $v=uw$ for some 
$u\in Z(U_{q}(f(K,H)))$.
\end{thm}
{\bf Proof:} If $v=uw$ for some $u \in Z(U_{q}(f(K,H)))$, 
then it is obvious that $v$ is a Whittaker vector. 
Conversely, let $v=uw$ for some $u\in U_{q}(f(K,H))$ 
be a Whittaker vector of $V$. Then $v=u^{\eta}w$ by 
the definition of Whittaker module. So we may 
assume that $u\in U_{q}(F,K^{\pm 1},H^{\pm 1})$. If $x\in U_{q}(E)$, 
then we have $xuw=\eta(x)uw$ and $uxw=\eta(x)uw$. 
Thus $[x,u]w=0$ and hence $[x,u]^{\eta}w=0$. But 
we have $x\bullet u=[x,u]^{\eta}$. Thus we have 
$u\in X$. We can now write $u=u_{1}+u_{2}$ with 
$u_{1}\in U_{q,w}(F,K^{\pm 1},H^{\pm 1})$ and $u_{2}\in W(F,K^{\pm
  1},H^{\pm 1})$. 
Then $u_{1}w=0$. Hence $u_{2}w=v$. But $u_{2}=u_{3}^{\eta}$ 
with $u_{3}\in Z(U_{q}(f(K,H)))$. So we have 
$v=u_{3}w$ which proves the theorem.
\qed
 
Now let $V$ be a Whittaker module and 
$End_{U_{q}(f(K,H))}(V)$ be the endomorphism 
ring of $V$ as a $U_{q}(f(K,H))-$module. 
Then we can define the following homomorphism 
of algebras using the action of $Z(U_{q}(f(K,H)))$ 
on $V$:
\[
\pi_{V}\colon Z(U_{q}(f(K,H))\longrightarrow End_{U_{q}(f(K,H))}(V).
\]
It is clear that 
\[
Z(U_{q}(f(K,H)))/Z_{V}(U_{q}(f(K,H)))\cong \pi_{V}(Z(U_{q}(f(K,H))))\subset End_{U_{q}(f(K,H))}(V).
\]
In fact, the next theorem says that this inclusion is an equality as well.
\begin{thm}
Let $V$ be a Whittaker $U_{q}(f(K,H))-$module. Then
$End_{U_{q}(f(K,H))}(V)\cong Z(U_{q}(f(K,H)))/Z_{V}$. 
In particular, $End_{U_{q}(f(K,H))}(V)$ is commutative.
\end{thm}
{\bf Proof:} Let $w\in V$ be a cyclic Whittaker vector. 
If $\alpha \in End_{U_{q}(f(K,H))}(V)$, then 
$\alpha(w)=uw$ for some $u\in Z(U_{q}(f(K,H)))$ by Theorem 4.5. 
Thus we have $\alpha(vw)=vuw=uvw$. Hence $\alpha=\pi_{u}$, which 
proves the theorem.
\qed

Now we are going to construct explicitly some Whittaker 
modules. Let 
\[
\xi \colon Z(U_{q}(f(K,H))) \longrightarrow \C
\] 
be a central character of the center $Z(U_{q}(f(K,H)))$. For any given
central character $\xi$, let $Z_{\xi}=Ker(\xi)\subset Z(U_{q}(f(K,H)))$ 
and $Z_{\xi}$ is a maximal ideal of $Z(U_{q}(f(K,H)))$. 
Since $\C$ is algebraically closed, then $Z_{\xi}=(\Omega
-a_{\xi},KH-b_{\xi})$ for some $a_{\xi} \in \C, b_{\xi}\in
\C^{\ast}$. For any given central character $\xi$, let $\C_{\xi,\eta}$
be the one dimensional $Z(U_{q}(f(K,H)))\otimes U_{q}(E)-$module
defined by $uvy=\xi(u)\eta(v)y$ for any $u\in Z(U_{q}(f(K,H)))$ and 
any $v\in U_{q}(E)$. We set
\[
Y_{\xi,\eta}=U_{q}(f(K,H))\otimes_{Z(U_{q}(f(K,H)))\otimes U_{q}(E)} \C_{\xi,\eta}.
\] 
It is easy to see that $Y_{\xi, \eta}$ is a Whittaker module of type
$\eta$ and admits a central character $\xi$. By Schur's lemma, we know 
every irreducible representation has a central character. As studied 
in \cite{JWS}, we know $U_{q}(f(K,H))$ has a similar theory for Verma 
modules. In fact, Verma modules also fall into the category of
Whittaker modules if we take the trivial character of
$U_{q}(E)$. Namely we have the following 
\[
M_{\lambda}=U_{q}(f(K,H))\otimes_{U_{q}(E,K^{\pm 1},H^{\pm 1})}\C_{\lambda},
\]
where $K, H$ act on $\C_{\lambda}$ through the character $\lambda$ of 
$\C[K^{\pm 1},H^{\pm 1}]$ and $U_{q}(E)$ act trivially on
$\C_{\lambda}$. Thus, $M_{\lambda}$ admits a central character. 
It is well-known that Verma modules may not be necessarily
irreducible, even though they have central characters. However, 
Whittaker modules are in the other extreme as shown in the next 
theorem:
\begin{thm} 
Let $V$ be a Whittaker module for $U_{q}(f(K,H))$. Then the following 
statements are equivalent.
\begin{enumerate}
\item $V$ is irreducible.\\
\item $V$ admits a central character.\\
\item $Z_{V}$ is a maximal ideal.\\
\item The space of Whittaker vectors of $V$ is one-dimensional.\\
\item All nonzero Whittaker vectors of $V$ are cyclic.\\
\item The centralizer $End_{U_{q}(f(K,H))}(V)$ is reduced to  $\C$.\\
\item $V$ is isomorphic to $Y_{\xi,\eta}$ for some central character $\xi$.
\end{enumerate}
\end{thm}
{\bf Proof:} It is easy to see that $(2)-(7)$ are equivalent to each
other by using the previous Theorems we have just proved. Since $\C$
is algebraically closed and uncountable, we also know $(1)$ implies
$(2)$ by using a theorem due to Dixmier \cite{Di}. To complete the 
proof, it suffices to show that $(2)$ implies $(1)$, namely if $V$ has
a central character, then $V$ is irreducible. 

Let $\omega \in V$ be a cyclic Whittaker vector, then
$V=U_{q}(f(K,H))\omega$. We have $V=U_{q}(F,K^{\pm 1},H^{\pm 1})w$. 
Since $V$ is irreducible, then $V$ has a central character. Thus 
we have $\Omega w=\lambda(\Omega)w$. Now we have 
\[
\Omega w=(\eta(E)F+\frac{q^{2m}K^{m}+H^{m}}{(q^{2m}-1)(q-q^{-1})})w.
\]
Hence the action of $F$ on $V$ is uniquely determined by the 
action of $K$ and $H$ on $V$, and $H^{-1}v=aKv, K^{-1}v=bHv$ 
for some $a, b\in \C^{\ast}$ and for any $v\in V$. Thus $V$ 
has a $\C-$basis consisting of elements 
$\{K^{i}\omega, H^{j}\omega \mid i,j \in \Z_{\geq 0}\}$. 

Let 
\[
0\neq v=(\sum_{i=0}^{n}a_{i}K^{i}+\sum_{j=1}^{m}b_{j}H^{j})\omega \in
V,
\] 
then
\begin{eqnarray*}
E(\sum_{i=0}^{n}a_{i}K^{i}+\sum_{j=1}^{m}b_{j}H^{j})\omega&=&(\sum_{i=0}^{n} 
q^{-2i}a_{i}K^{i}+\sum_{j=1}^{m}q^{2j}b_{j}H^{j})E\omega\\
&=&\eta(E)(\sum_{i=0}^{n}q^{-2i}a_{i}K^{i}+\sum_{j=1}^{m}q^{2j}b_{j}H^{j} )\omega.
\end{eqnarray*}

Thus we have $0 \neq \eta(E) q^{-2n}v-Ev\in V$, in which the top degree of $K$ 
is $n-1$. By repeating this operation finitely many times, we will 
finally get an element $0\neq a\omega$ with $a\in \C^{\ast}$. This 
means that $V=U_{q}(f(K,H))v$ for any $0\neq v\in V$. So $V$ is
irreducible. Therefore, we are done with the proof. \qed

In addition, the proof of the previous theorem also implies the following:
\begin{thm}
Let $(V,w)$ be an irreducible Whittaker module 
with a Whittaker vector $w$, then $V$ has a 
$\C-$basis consisting of elements $\{K^{i}\omega,\, H^{j}\omega\mid i,j \in
\Z_{\geq 0} \}$.
\end{thm}
 \qed

It is easy to show the following two theorems, for more details about
the proof, we refer the reader to \cite{K}.
\begin{thm}
Let $V$ be a $U_{q}(f(K,H))-$module which admits a central
 character. Assume that $w\in V$ is a Whittaker vector.
Then the submodule $U_{q}(f(K,H))w\subset V$ is irreducible.
\end{thm}
\qed

\begin{thm}
Let $V_{1},V_{2}$ be any two irreducible $U_{q}(f(K,H))-$modules 
with the same central character. If $V_{1}$ and $V_{2}$ 
contain Whittaker vectors, then these vectors are unique 
up to scalars. And furthermore, $V_{1}$ and $V_{2}$ are 
isomorphic to each other as $U_{q}(f(K,H))-$modules.
\end{thm}
\qed

\end{document}